\documentclass[reqno,12pt]{amsart} %
\usepackage{amsmath,amstext, amsthm, amssymb,amsfonts}
\usepackage{fancybox}
\numberwithin{equation}{section}

\def\RR{{\mathbb R}}

\def\NN{{\mathbb N}}

\def\eps{\varepsilon}
\newcommand{\la}{\lambda}

\newcommand{\calL}{\mathcal{L}}
\newcommand{\calF}{\mathcal{F}}

\newcommand{\calD}{\mathcal{D}}

\newtheorem{theorem}{Theorem}[section]
\newtheorem{proposition}[theorem]{Proposition}

\newtheorem{corollary}[theorem]{Corollary}
\theoremstyle{definition}
\newtheorem{definition}{Definition}[section]
\newtheorem{remark}{Remark}[section]
\newtheorem{example}{Example}[section]

\newtheorem{conjecture}{Conjecture}[section]

\newcount\mins \newcount\hours \hours=\time \mins=\time
\divide\hours by60 \multiply\hours by 60 \advance\mins by-\hours
\divide\hours by60

\def\now{ \ifnum\hours>11 \ifnum\hours>12 \advance\hours by
-12 \fi \number\hours:\ifnum\mins<10 0\fi \number\mins\ pm,\ \else
\ifnum\hours=0 \hours=12 \fi \number\hours:\ifnum\mins<10 0\fi
\number\mins\ am,\ \fi}

\newcommand{\ToDoList}[1]{\shadowbox{\begin{minipage}{3.5in} #1 \end{minipage}}}
\newcommand{\Comment}[1]{\marginpar{\footnotesize #1}}
\renewcommand{\ToDoList}[1]{}
\renewcommand{\Comment}[1]{}
\title{Free  Exponential Families as Kernel Families}

\author{
W{\l}odzimierz  Bryc}\thanks{\noindent Research partially supported by NSF
grant  \#DMS-05-04198 and by Taft Research Seminar 2008/09}

\date{ September 26, 2008. Revised: April 24, 2009.
}

\address{
Department of Mathematical Sciences,\\
University of Cincinnati,\\
PO Box 210025,\\
Cincinnati, OH 45221--0025, USA}
\email{Wlodzimierz.Bryc@UC.edu}
\keywords{exponential families, dispersion models, free convolution, free cumulants}

\subjclass[2000]{60E10; 46L54}
\begin{document}


\begin{abstract}
Free exponential families have been previously introduced  as a
 special case of the $q$-exponential family.
We show that free exponential families arise also from the approach
analogous to the definition of exponential families by using  the Cauchy-Stieltjes kernel $1/(1-\theta
x)$ instead of the exponential kernel $\exp({\theta x})$. We use this approach to re-derive some known results   and
to study further similarities with exponential
families and reproductive exponential models.
\end{abstract}
\maketitle

\section{Introduction}

Since the seminal work of Voiculescu \cite{Voiculescu86}, there has
 been a flurry of activity on how
properties of free convolution
$\mu\boxplus\nu$ of probability measures are similar to and how they differ from  properties
of classical convolution $\mu*\nu$.
In particular, free probability analogues of the Central Limit Theorem, of the Poisson limit theorem, and
the L\'evy-Khinchin representation of $\boxplus$-infinitely divisible laws are now known, see \cite{Hiai-Petz00}.
New additional analogies between free and classical probability are developed in \cite{Ben_Arous:2009,Ben-Arous:2006}. In this paper we study a free probability analogue of the concept of exponential family.

 Free exponential families were introduced in \cite[Definition 4.1]{Bryc-Ismail-05} as part of a study of the relations between approximation
 operators, classical exponential families and their $q$-deformations.
An alternative approach to free exponential families  which  we adopt in this paper
emphasizes similarities to classical exponential families, and is based on an idea of kernel family
 introduced in \cite{Wesolowski90}.
 We show that the two approaches are closely related, and
that  every non-degenerate
compactly supported probability measure generates a free exponential family, see  Theorem \ref{T1}.
We then  relate  variance functions of free exponential families to  free cumulants.
This relation is simpler than the  corresponding relation for classical
 exponential families and is expressed by a concise formula.
We apply the formula  to compute  free cumulants of the "free gamma" law which were
 stated without proof in
\cite{Bozejko-Bryc-04},   to derive  simple necessary conditions for a smooth  function to be the variance function of a free exponential family,   and to  investigate similarities with
classical dispersion models \cite{Jorg}.

\section{Cauchy-Stieltjes Kernel Families}\label{Kernel Section}
According to Weso\l owski \cite{Wesolowski90}, the kernel family generated by a
kernel $k(x,\theta)$ consists of the probability measures
$$\{k(x,\theta)/L(\theta)\nu(dx):\; \theta\in\Theta\},$$ where
$L(\theta)=\int k(x,\theta)\nu(dx)$ is the normalizing constant, and
$\nu$ is the generating measure.

The theory of exponential families is based on
the kernel $k(x,\theta)=e^{\theta x}$. See, e.g., \cite{Jorg}, \cite{Let},  or \cite[Section 2.3]{Diaconis-08}.
In this paper, we consider the Cauchy-Stieltjes  kernel
\begin{equation*}
k(x,\theta)=\frac{1}{1-  \theta x}.
\end{equation*}

\begin{definition}
  Suppose $\nu$ is a compactly supported non-degenerate (i.e. not a point mass) probability measure. Let \begin{equation*}  
M(\theta)=\int \frac{1}{1- \theta x} \nu(dx).
\end{equation*} The Cauchy-Stieltjes family generated by
$\nu$ is the family of probability measures
\begin{equation}\label{Kernel family}
\mathcal{K}(\nu;\Theta)=\left\{P_\theta(dx)=\frac{1}{M(\theta)(1-\theta
x)}\nu(dx): \theta\in\Theta\right\},
\end{equation}
where $\Theta\ni 0$ is an open set on which $M(\theta)$ is well
defined, strictly positive and $\theta\; \rm{supp}(\nu)\subset(-\infty,1)$. (We shall only consider $\Theta=(-\eps,\eps)$ with $\eps>0$ small
enough.)
\end{definition}
Our first goal is to show that the Cauchy-Stieltjes family is essentially the same concept as the concept of free exponential family introduced in \cite{Bryc-Ismail-05}.
We begin with  a suitable reparametrization of $\mathcal{K}(\nu;\Theta)$.
\subsection{Parameterizations by the mean}
From \eqref{Kernel family} we compute the mean
$m(\theta)=\int x P_\theta(dx)$. Since $P_0=\nu$ we get $m(0)=\int
x\nu(dx)=m_0$ and for $\theta\ne 0$, a calculation gives
\begin{equation}\label{L2m}
m(\theta)=\frac{M(\theta)-1}{\theta M(\theta)}.
\end{equation}
Since $M(0)=1$ and $M(\theta)$ is analytic at $\theta=0$,  we see
that $m(\theta)$ is analytic for $|\theta|$ small enough. We have
$$m'(\theta)=\frac{M(\theta)+\theta M'(\theta)-(M(\theta))^2}{\theta^2(M(\theta))^2}.$$
Since $\nu$ is non-degenerate,
\begin{multline}\label{increasing}
M(\theta)+\theta M'(\theta)-(M(\theta))^2\\
=\int\frac{1}{(1-\theta x)^2}\nu(dx)-\left(\int\frac{1}{1-\theta x}\nu(dx)\right)^2>0
\end{multline}
for all $|\theta|>0$ small enough. Thus the function $\theta\mapsto
m(\theta)$ is increasing on an open interval containing $0$. Denoting by
$\psi$  the inverse function, we are thus lead to
parametrization of a subset of $\mathcal{K}(\nu,\Theta)$  by the mean,
\begin{equation}\label{Kernel family2}
\mathcal{L}(\nu;R)=\left\{Q_m(dx)=P_{\psi(m)}(dx): m\in R\right\},
\end{equation}
where $R=m(\Theta_0)$, and $\Theta_0\subset\Theta$ is an
appropriate interval with $0\in\Theta_0$. Notice that we refrain from
claiming that \eqref{Kernel family} and \eqref{Kernel family2} are
equivalent: we only claim that for any pair of open sets $\Theta\ni0$ and $R\ni m_0$, there are open  sets  $\Theta_0\ni 0$ and $ R_0\ni m_0$
such that
$\mathcal{L}(\nu;R_0)\subset \mathcal{K}(\nu;\Theta)$,  { and }
$\mathcal{K}(\nu;\Theta_0)\subset \mathcal{L}(\nu;R)$.

The variance function of the Cauchy-Stieltjes family \eqref{Kernel
family2} is
\begin{equation}
  \label{Def Var}
  V(m)=\int (x-m)^2 Q_m(dx).
\end{equation}


\section{Relation to free exponential families}
 The following generalizes slightly \cite[Section 4]{Bryc-Ismail-05}; note that this definition is not constructive: for a given $V$, the corresponding  free exponential  family may fail to exist, see Example \ref{Ex_no_FEF}.
 \begin{definition}
 The free exponential family with variance function $V$ generated by a  compactly supported measure $\nu$ with mean $m_0\in(A,B)$ is a family of probability measures
\begin{equation}
  \label{F(V)}
\calF_{m_0}(V)=\left\{\frac{V(m)}{V(m)+(m-m_0)(m-x)}\nu(dx):\;
m\in(A,B)\right\}.
\end{equation}
\end{definition}

The next result shows that Cauchy-Stieltjes kernel families under
parametrization by the mean are essentially the same as free exponential
families, thus providing existence argument for  free exponential families. Furthermore, the generating
 measure $\nu$ is determined uniquely by $m_0$ and the variance function $V(m)$;
 the latter is an analog of
the classical uniqueness theorem  for exponential families, see
\cite[Theorem 2.11]{Jorg} or \cite[Proposition 2.2]{Let:Mor}.

Recall that the Cauchy-Stieltjes transform of a probability measure
$\nu$ is
\begin{equation}\label{G-transform}
G_\nu(z)=\int\frac{1}{z-x}\nu(dx).
\end{equation}
If $\nu$ is compactly supported then $G_\nu$ is analytic in the
neighborhood of $\infty$ in the complex plane; in particular,
compactly supported measures are determined uniquely by $G_\nu(z)$
for large enough real $z$.
\begin{theorem}
  \label{T1}
  Let $\{P_\theta: \theta\in\Theta\}$ be a Cauchy-Stieltjes family \eqref{Kernel family} generated by a non-degenerate compactly supported measure $\nu$ with $\int x \nu(dx)=m_0$. Then there is  a neighborhood $(A,B)$ of $m_0$
  in which the variance function $V$ in \eqref{Def Var}
  is analytic, strictly positive, and
  \begin{equation}
P_{\psi(m)}(dx) = \frac{V(m)}{V(m)+(m-m_0)(m-x)}\nu(dx),\; m\in(A,B).
\end{equation}

 Conversely, if $V$ is analytic and strictly
  positive  in a neighborhood of $m_0$, and there is a probability measure $\nu$ with mean $m_0$ such that the (positive) measures
     \begin{equation}
  \label{Q_m}
Q_m(dx)=\frac{V(m)}{V(m)+(m-m_0)(m-x)}\nu(dx)
\end{equation}
are probability measures for all $m$ in a neighborhood of $m_0$, then $\nu$ is compactly supported, non-degenerate, and
is  determined uniquely by \eqref{G-transform} with
  \begin{equation}
  \label{G2V}
  G_\nu(z)=\frac{m-m_0}{V(m)},
\end{equation}
with $z=1/\psi(m)=m+V(m)/(m-m_0)$. In particular, in a perhaps smaller neighborhood $R$ of $m_0$, probability measures \eqref{Q_m}
for $ m\in R$
are included in the  Cauchy-Stieltjes family \eqref{Kernel
  family2} generated by $\nu$.
\end{theorem}
\begin{proof}
 We first calculate the variance $v(\theta)=\int x^2
P_\theta(dx)-m^2(\theta)=\int x(x-m_0)P_\theta(dx)-m(\theta)(m(\theta)-m_0)$. Since
\begin{multline*}
\int x(x-m_0) P_\theta(dx)=\int
\frac{(x-m_0)(\theta x-1)+x-m_0}{\theta M(\theta)(1-\theta x)}\nu(dx)
\\=\frac{m(\theta)-m_0}{\theta},
\end{multline*}
we see that for $\theta \ne 0$  the variance is
\begin{equation}\label{m2v}
v(\theta)=(m(\theta)-m_0)\left(\frac{1}{\theta}-m(\theta)\right).
\end{equation}
Since $m(\theta)$ is analytic at $\theta=0$ and $m(0)=m_0$, this shows that
that $v(\theta)$ is analytic at $\theta=0$.
 Let $V(m)=v(\psi(m))$ denote the variance function in
parametrization of (a subset of) $\mathcal{K}$ by the mean; clearly
$V$ is an analytic function in a neighborhood of $m=m_0$.

Now  \eqref{m2v} implies $(m(\theta)-m_0)/v(\theta)=\theta/(1-\theta m(\theta))$ and \eqref{L2m} implies $1-\theta m(\theta)=1/M(\theta)$,
so \eqref{Q_m} is
equivalent to
\begin{multline*}
Q_{m(\theta)}(dx)=\frac{1}{1+(m(\theta)-x)(m(\theta)-m_0)/v(\theta)}\nu(dx)
\\
=\frac{1-\theta m(\theta)}{1-\theta m(\theta)+(m(\theta)-x)\theta}\nu(dx)
=\frac{1-\theta m(\theta)}{1-\theta x}\nu(dx)\\
=\frac{1}{M(\theta)(1-\theta x)}\nu(dx)=P_\theta(dx).
\end{multline*}


To prove the converse implication, note that for $m$ such that
$V(m)>0$  we can re-write
 $\int Q_m(dx)=1$ as
 $$\int \frac{1}{V(m)/(m-m_0)+m-x}\nu(dx)= \frac{m-m_0}{V(m)}.$$
 Thus with
 \begin{equation}\label{z2m}
 z=m+\frac{V(m)}{ m-m_0},
 \end{equation} we get
 \eqref{G2V}. Since $\lim_{m\to m_0^\pm}V(m)/(m-m_0)=\pm \infty$, this shows that Cauchy-Stieltjes transform $G_\nu(z)$ is defined for all real $z$ with $|z|$ large enough. This implies that
 $\nu$ has compact support, with moments that are uniquely determined from the corresponding moment generating function $M(z)=1/z G(1/z)$ for $z$ small enough.  (Compactness of support is also proved  more directly in the proof of Theorem \ref{V2c}.) Finally, $\nu$ is non-degenerate as its variance is $V(m_0)>0$.

\end{proof}

\begin{remark}\label{Rem-not-V}
  Solving equation \eqref{m2v} for $\theta$ we see that $$\psi(m)=\frac{m-m_0}{m(m-m_0)+V(m)}.$$ Thus
  a necessary condition for $V$ to be a variance function is that $m\mapsto m+V(m)/(m-m_0)$ is decreasing in a neighborhood of $m_0$, see \eqref{increasing}.
\end{remark}
 \subsection{Free exponential families with quadratic variance function}
In this section we recall  \cite[Theorem 4.2
]{Bryc-Ismail-05}; since manuscript  \cite{Bryc-Ismail-05} is available  in preprint form only and we have already set up all identities needed for the proof, we include the argument which is taken   from \cite{Bryc-Ismail-05}. The corresponding result for classical
exponential families is \cite[Theorem 3.3]{Ism:May} and
\cite[Section 4]{Mor}; the  result for $q$-exponential families is \cite[Theorem 3.2
]{Bryc-Ismail-05}.
\begin{theorem}[{\cite[Theorem 4.2
]{Bryc-Ismail-05}}]\label{T free quadr}
Suppose $b\geq -1$,  $m_0=0$. The free exponential family with the
variance function
$$V(m)=1+a m+b m^2$$ consists of probability measures \eqref{Q_m}
with generating measure  
\begin{multline}\label{free mu}
\nu(dx)=\frac{\sqrt{4(1+b)-(x-a)^2}}{2\pi (bx^2 +ax +1)}
1_{(a-2\sqrt{1+b}, a+2\sqrt{1+b})} dx +p_1\delta_{x_1}+
p_2\delta_{x_2},
\end{multline}
where the discrete part of $\nu$ is absent except for the following cases:
\begin{enumerate}
\item if $b=0, a^2>1$, then
$ p_1=1-1/a^2$, $x_1={-1/a}$, $p_2=0$.
\item if $b>0$ and $a^2>4 b$, then
$ p_1=\max\left\{0, 1-\frac{|a|-\sqrt{a^2-4b}}{2
b\sqrt{a^2-4b}}\right\} $, $p_2=0$, and $x_1=
\pm\frac{|a|-\sqrt{a^2-4b}}{2b}$ with the sign opposite to the sign
of $a$.
\item  if $-1\leq b<0$ then there are two atoms at %
$$
x_{1,2}=\frac{-a\pm \sqrt{a^2-4b}}{2b}, \; p_{1,2}=
1+\frac{\sqrt{a^2-4b}\mp a}{2 b\sqrt{a^2-4b}}.
$$
\end{enumerate}
\end{theorem}
\begin{proof}
With  $m_0=0$  and
$V(m)=1+a m + b m^2$, equation \eqref{z2m} can be solved for $m$,
giving
$$
m=\frac{z-a - \sqrt{\left( a - z \right)^2-4\,\left( 1 + b
\right)}}{2\,
         \left( 1 + b \right) } ,
             $$
             so \eqref{G2V} gives
\begin{equation}\label{Free meixner G}
G_\nu(z)= \frac{a + z + 2\,b\,z- {\sqrt{
         {\left( a - z \right) }^2-4\,\left( 1 + b \right)  }} }{2\,
     \left( 1 + a z\ + b\,z^2  \right) }.
\end{equation}
This Cauchy-Stieltjes transform corresponds to the free-Meixner law \eqref{free mu}, see
\cite{Anshelevich01,Saitoh-Yoshida01}.
\end{proof}

Theorem \ref{T free quadr} results covers a number of important laws that appeared in the literature.
Up to a dilation and
convolution with a degenerate law $\delta_a$ (i.e. up to "the type")
the generating measure $\nu$ is:
\begin{enumerate}
\item  the Wigner's semicircle (free Gaussian) law if $a=b=0$; see \cite[Section 2.5]{Voiculescu00};
\item  the Marchenko-Pastur (free Poisson) type law if $b=0$ and $a\ne 0$; see \cite[Section 2.7]{Voiculescu00};
\item the free Pascal (free negative binomial) type law if $b>0$ and $a^2>4b$; see \cite[Example 3.6]{Saitoh-Yoshida01};
\item the free Gamma type law if $b>0$ and $a^2=4b$;  see \cite[Proposition 3.6]{Bozejko-Bryc-04};
\item the free analog of hyperbolic type law if $b> 0$ and $a^2<4b$; see \cite[Theorem 4]{Anshelevich01};
\item the free binomial type law    if $-1\leq b <0$; see \cite[Example 3.4]{Saitoh-Yoshida01} and
\cite[Proposition 2.1]{Bozejko-Bryc-04}.
\end{enumerate}
  The laws in (i)-(v) are infinitely divisible with respect  to free additive convolution
 (we recall the definition near \eqref {free convolution}).
In \cite[Theorem 4]{Anshelevich01} they appear in connection to
martingale polynomials with respect to
  free L\'evy processes; free infinite divisibility is analyzed also in  \cite{Saitoh-Yoshida01};  \cite{Anshelevich04} studies further free probability aspects of this family;
  in  \cite[Theorem 3.2]{Bozejko-Bryc-04}
the same laws appear as a solution  to a  quadratic regression
problem in free probability;
 in
  \cite[Theorem 4.3]{Bryc-Wesolowski-03} these laws occur in a ``classical regression"
  problem.

\subsection{Free Cumulants and Variance Functions}\label{Sect FC2VF}

Recall that if $\nu$ is a compactly supported measure with the
Cauchy-Stieltjes transform $G_\nu$, then the inverse function
$K_\nu(z)=G^{-1}_\nu(z)$  exists for small enough $z\ne 0$, see
\cite{Voiculescu00}. The $R$-transform is defined as
\begin{equation}\label{R def} R_\nu(z)=K_\nu(z)-1/z \end{equation}
and is analytic at $z=0$, \begin{equation}\label{R}
R_\nu(z)=\sum_{n=1}^\infty c_n(\nu) z^{n-1}. \end{equation} The
coefficients $c_n=c_n(\nu)$ are called free cumulants of measure $\nu$, see
 \cite{Speicher-97}.

  The following  result extends {\cite[Remark 4.4]{Bryc-Ismail-05}} and plays a  role  analogous to J\o rgensen's theorem \cite[Theorem 3.2]{Let:Mor}. (For the formula connecting classical cumulants with
the variance functions of natural exponential families, see
\cite[(2.10)]{Mor} or \cite[Exercise 2.14]{Jorg}.)
\begin{theorem}\label{V2c}\tolerance=1000
Suppose $V$ is analytic in a neighborhood of $m_0$, $V(m_0)>0$, and
$\nu$ is a probability measure with finite all moments, such that
$\int x \nu(dx)=m_0$.
 Then the following conditions are equivalent.
\begin{enumerate}
\item  $\nu$ is  non-degenerate, compactly supported, and there exists an interval $(A,B)\ni m_0$ such that
 \eqref{F(V)} defines a family of probability measures parameterized by the mean with  the variance
 function $V$.
 \item \label{V2cii} The free cumulants \eqref{R} of $\nu$ are
$c_1=m_0$,  and for $n\geq 1$
\begin{equation}\label{free cumulants}
c_{n+1}=\left.\frac{1}{ n!}\frac{ d^{n-1}}{d
x^{n-1}}\left(V(x)\right)^n\right|_{x=m_0}.
\end{equation}

\end{enumerate}
\end{theorem}

\begin{proof} 
 Suppose  that $V$ determines the free exponential family generated
by a compactly supported measure $\nu$.  For $m\ne 0$ close
enough to $0$, from \eqref{G2V} and \eqref{z2m} we get
$$
m+\frac{V(m)}{m-m_0}=K_\nu\left(\frac{m-m_0}{V(m)}\right).
$$
Thus,  \eqref{R def} says that the $R$-transform of $\nu$ satisfies
\begin{equation}\label{R formula}
R_\nu\left(\frac{m-m_0}{V(m)}\right)=m.
\end{equation}
From this we derive \eqref{free cumulants} by using  the Lagrange
expansion theorem, which says that if $\phi(z)$ is analytic
in a neighborhood of $z=m_0$, $\phi(m_0)\ne 0$ and $\xi :=
{(m-m_0)}/{\phi(m)}$ then
\begin{equation} \label{eqlagrange} m(\xi) = m_0 +
\sum_{n=1}^\infty \frac{\xi^n}{n!} \, \left[\frac{d^{n-1}
[\phi(x)]^n}{dx^{n-1}}\right]_{x=m_0}. \end{equation} (See, e.g.,
\cite[(L), page 145]{Pol:Sze}.)

{\tolerance=500
Suppose now that a probability measure $\nu$ satisfies \eqref{free
cumulants} and $\int x\nu(dx)=m_0$. Then  the variance  $c_2(\nu)=V(m_0)>0$,  so $\nu$  it is non-degenerate.  We first verify that $\nu$ has
compact support. Since $V$ is analytic, \eqref{free cumulants} implies that
\begin{equation}
\label{Cauchy}
c_{n+1}=\frac1{2\pi n i }\oint_{|z-m_0|=\delta} V(z)^n/(z-m_0)^{n-1},
\end{equation}
 so  there exist $M>0$ such
that $|c_n|\leq M^n$ for all $n\geq 1$. The compactness of support follows now from
\cite[Corollary 1.6]{Benaych-Georges04}; for  completeness we include the proof.
Denoting by $\mathcal{NC}[n]$ the set of non-crossing partitions of $\{1,2,\dots,n\}$,
from   \cite[(2.5.8)]{Hiai-Petz00} we have
$$
\int x^{2n}\nu(dx)=\sum_{\mathcal{V}\in \mathcal{NC}[2n]}
\prod_{B\in\mathcal{V}} c_{|B|}\leq M^{2n}\# \mathcal{NC}[2n]=M^{2n}\frac{1}{2n+1}
\binom{4n}{2n};
$$
for the last equality, see \cite[(2.5.11)]{Hiai-Petz00}. Since the $m$-th Catalan number is less than $4^m$,
 $$\limsup_{n\to\infty}\left(\int |x|^{2n}\nu(dx)\right)^{1/(2n)}\leq 4M<\infty,$$
and $\nu$ has compact support.}

From $\mbox{supp}(\nu)\subset[-4M,4M]$ we deduce that the  Cauchy-Stieltjes transform $G_\nu(z)$ is analytic for $|z|>4M$, and
the $R$-series is analytic for all $|z|$ small enough.

Since $V(m)>0$ for $m$ close enough to $m_0$,  taking the derivative we see that $z\mapsto {(z-m_0)}/{V(z)}$ is increasing in
a neighborhood of $z=m_0$. Denoting by $h$ the inverse, we have
$$
h\left(\frac{z-m_0}{V(z)}\right)=z.
$$
From $c_1(\nu)=m_0$ we see that
 $R(0)=m_0=h(0)$. By \eqref{eqlagrange},
  we see that all derivatives of $h$ at $z=0$ match the derivatives of $R$.
 Thus $h(z)=R(z)$ and
\eqref{R formula} holds for all $m$ in a neighborhood of $0$. For
analytic $G_\nu$, the latter is equivalent to \eqref{Q_m} holding
for all $m$ close enough to $0$. Thus $V(m)$ is the variance
function of a
free exponential family generated by $\nu$ with  $m\in(-\delta,\delta)$ for some $\delta>0$.
\end{proof}


We  now  use \eqref{free cumulants} to relate certain free cumulants to Catalan numbers.
\begin{corollary} \cite[Remark 5.7]{Bozejko-Bryc-04}
 If $\nu$ is the standardized  free gamma Meixner law, i.e. it generates the free exponential family
 with $m_0=0$ and  variance function $V(m)=(1+a m)^2$, then its free cumulants are $$c_{k+1}(\nu)=\frac{1}{k+1}\binom{2k}{k} a^{k-1},\;k\geq 1.$$
\end{corollary}
This fact was stated without proof in \cite[Remark 5.7]{Bozejko-Bryc-04}; the approach indicated there lead to a relatively long proof.
\begin{proof}
From \eqref{free cumulants},
\begin{multline*}
c_{k+1}(\nu)=\left.\frac{1}{  k!}\frac{d^{k-1}}{d
x^{k-1}}\left(1+ax\right)^{2k}\right|_{x=0} =a^{k-1}\frac{2k
(2k-1)\dots (k+2)}{k!}\\=\frac{a^{k-1}}{k+1}\frac{(2k)!}{(k!)^2}\;.
\end{multline*}

\end{proof}

Recall that the free additive convolution of compactly supported probability measures $\mu,\nu$ is a unique
compactly supported measure denoted by $\mu\boxplus\nu$   with the
$R$-transform
\begin{equation}\label{R-sum}
R_{\mu\boxplus\nu}(z)=R_\mu(z)+R_\nu(z).
\end{equation}
(See \cite{Voiculescu86}.)
Equivalently, free cumulants linearize free convolution,
\begin{equation}\label{free convolution}
c_n(\mu\boxplus \nu)=c_n(\mu)+c_n(\nu),\; n\geq 1
\end{equation}
 just like classical cumulants linearize the classical convolution.
 Recall that $\mu$ is $\boxplus$-infinitely divisible if for every $n=1,2,\dots,$ there is a measure $\nu$ such that $\mu=\nu\boxplus\nu\boxplus\dots\boxplus\nu$ (the $n$-fold free convolution).

\begin{corollary}\label{C.4.5}   $V(m)=1/(1-m)$ is a variance function of a free exponential family
generated by the centered $\boxplus$-infinitely divisible   measure $\nu$ with free cumulants
 $$c_{k+1}(\nu)=\frac{1}{k}\binom{2k-2}{k-1}, \;k\geq 1.$$
\end{corollary}
\begin{proof}
From \eqref{free cumulants},
$$c_{k+1}=\left.\frac{1}{  k!}\frac{d^{k-1}}{d x^{k-1}}\left(1-x\right)^{-k}\right|_{x=0}=
\frac{k(k+1)\dots(2k-2)}{k!}=
\frac{1}{k} \frac{(2k-2)!}{((k-1)!)^2}\;.
$$
It is well known that Catalan numbers are even moments of the semicircle law,
 $$c_{k+1}=\int_{-2}^{2}x^{2k}\frac{\sqrt{4-x^2}}{2\pi}dx=\int_{0}^{{4}}x^{k}\frac{\sqrt{4/x-1}}{2\pi}dx\;.$$
Therefore, $R(z)=\sum_{k=1}^\infty c_k
z^{k-1}=\int_{0}^{{4}}\frac{zx}{1-zx}\frac{\sqrt{4/x-1}}{2\pi}dx=
\int_{0}^{{4}}\frac{z}{1-zx}\frac{\sqrt{4x-x^2}}{2\pi}dx$
corresponds to $\boxplus$-infinitely divisible law, see \cite[Theorem
3.3.6]{Hiai-Petz00}. Thus Catalan numbers $c_{k+1}$ with $c_1=0$ are
indeed free cumulants of some  $\boxplus$-infinitely divisible  measure $\nu$.
\end{proof}

It is known that not every function $V$ is a variance function of a
natural  exponential family. It is therefore not surprising that not
every analytic functions $V$ can serve as the variance functions for
a free exponential family.
\begin{corollary}\label{C.4.6} Suppose $V$ is analytic at $0$ and $V(0)=1$,  $V''(0)<-2$.
Then  $V$ cannot be a variance function of a free exponential family
with $m_0=0$.
\end{corollary}
\begin{proof}
Suppose $V$ generates a free exponential  family with generating
measure $\nu$.  Let $m_j=\int x^j \nu(dx)$ with $m_1=0$, $m_2=1$.
Then the $3\times 3$  Hankel determinant is
$$
\det\left[ {\begin{array}{*{20}c}
   1 & m_1 & m_2  \\
   m_1 & m_2 & {m_3}  \\
   m_2 & {m_3} & {m_4}  \\
 \end{array} } \right] =\det\left[ {\begin{array}{*{20}c}
   1 & 0 & 1  \\
   0 & 1 & {m_3}  \\
   1 & {m_3} & {m_4}  \\
 \end{array} } \right] = m_4-m_3^2-1\geq 0.
$$
Using \eqref{free cumulants}, the fourth moment is
$$m_4=c_4(\nu)+2 c_2^2(\nu)=c_4(\nu)+2=V'(0)^2+\frac1{2}V''(0)+2$$ and
 $m_3=c_3(\nu)=V'(0)$, see \cite[(2.5.8)]{Hiai-Petz00}. Thus $m_4-m_3^2-1\geq 0$ translates into
$V''(0)\geq -2$.
\end{proof}

\begin{example}
 If $b<-1$, then $V(m)=1+am+bm^2$ is not a variance function of a free exponential  family with $m_0=0$.
 Compare Theorem \ref{T free quadr}.  (This can also be seen from Remark  \ref{Rem-not-V}.)
\end{example}\label{Ex_no_FEF}
\begin{example} $V(m)=(1-m)/(1+m)$ is not a variance function of a free exponential  family with $m_0=0$. (This can also be seen from Remark  \ref{Rem-not-V}.)
\end{example}

Combining \eqref{free cumulants} with the $\boxplus$-L\'evy-Khinchin
formula \cite[Theorem 3.3.6]{Hiai-Petz00}, compare \cite[Lemma
3.4]{Barndorff-Nielsen-Thorbj02a}, we get also the following.

\begin{corollary} Suppose $V(m)$ is analytic at $0$, $V(0)=1$. Then  the following conditions are equivalent.
\begin{enumerate}
\item There exists a centered $\boxplus$-infinitely divisible
probability measure $\nu$ such that $V$ is the variance function of
a free exponential  family generated by $\nu$.
\item There exists a  compactly supported probability measure $\omega$ such that
\begin{equation}\label{V2Levy}
\frac{1}{n!}\left.\frac{d^{n-1}}{d x^{n-1}}\left(V(x)\right)^n\right|_{x=0}=\int x^{n-1} \omega(dx), \;n\geq 1.
\end{equation}
\end{enumerate}
\end{corollary}
The Cauchy-Schwarz  inequality applied to the right hand side of the  L\'evy-Khinchin formula \eqref{V2Levy} implies
$(V^3)''/6\geq \left((V^2)'\right)^2/4$. This gives a simple necessary condition.
\begin{corollary}
If $V$ is analytic at $0$, $V(0)=1$, $V''(0)<0$ then $V$ cannot be
the variance function of  a free exponential  family generated by a
centered $\boxplus$-infinitely divisible measure.
\end{corollary}
We remark that the bound is sharp: from Theorem \ref{T free quadr}
we see that
 $V(m)=1$ is a variance function of the free exponential  family generated by
 the semicircle law; all of its members are  infinitely divisible, see Example \ref{Example: semicircle}.
\begin{example}[{Compare \cite[Theorem 3.2]{Saitoh-Yoshida01}}] If $b<0$ then $V(m)=1+am + bm^2$
cannot be the variance function of  a free exponential  family
generated by a centered $\boxplus$-infinitely divisible measure.
\end{example}

\subsection{Reproductive property}

Natural exponential families have two "reproductive" properties. The
first one is usually not named, and says that if a  compactly
supported  measure $\nu$ generates natural exponential family
$\calF$ and $\mu\in\calF(\nu)$ then $\calF(\mu)=\calF$. This is
usually interpreted as a statement that the natural exponential
family $\calF$ is determined solely
 by the variance function $V$ and can have many generating measures.

The analog of this property fails for free exponential families due
to the fact that the generating measure is determined uniquely
 by the variance function and parameter $m_0$. For example, a free exponential family $\calF$
generated by the centered semicircle law consists of the affine
transformations of the Marchenko-Pastur  laws, and for $m_0\ne 0$
the free exponential  family  generated by $\mu\in\calF$ with mean
$m_0$
 contains no other measures in common with $\calF$ except for $\mu$.

The second property which in \cite[(3.16)]{Jorg} is indeed called
the reproductive  property of an exponential family states that if
$\mu\in\calF(V)$, then for all $n\in\NN$ the law of the sample mean,
$D_n(\mu^{*n})$, is in $\calF(V/n)$.  Here  $D_r(\mu)(U):=\mu(r U)$
denotes the dilation of measure $\mu$ by a number $r\ne0$; in
probabilistic language, if $\mathcal{L}(X)=\mu$ then
$\mathcal{L}(X/r)=D_r(\mu)$.

Our goal is to prove an analogue of this result for the
Cauchy-Stieltjes families.

 Let $\mu^{\boxplus r}$  denote the $r$-fold free additive convolution of $\mu$ with itself. In contrast to classical convolution, this
 operation is well defined for all real
$r\geq 1$, see \cite{Nica-Speicher}.

\begin{proposition}[{\cite[Proposition 4.3]{Bryc-Ismail-05}}]\label{C.4.3} If a function $V$ analytic at $m_0$ is a variance function
of a free exponential family generated by a compactly supported
probability measure $\nu$ with    $m_0=\int x \nu(dx)$, then for
each $\la\geq 1$ there exists a neighborhood of $m_0$ such that
$V/\la$ is the variance function of the free exponential family
generated by measure
$$\nu_{\la}:=D_\la(\nu^{\boxplus \la}).$$
Moreover, if for each $\la>0$, there is a neighborhood of $m_0$ such
that  $V/\la$
 is a variance function of some free exponential family, then $\nu$ is $\boxplus$-infinitely divisible.
\end{proposition}

We note that in contrast to classical natural exponential families,
the neighborhood of $m_0$ where $m\mapsto V(m)/\la$ is a variance function
may vary
 with $\la$, see Example \ref{Example: semicircle}.

\ToDoList{How does the free sample averaging $D_n(\nu^{\boxplus n})$
change the interval $(A,B)$?
\begin{conjecture} If $\calF$ is defined over the interval
$(-\delta,\delta)$ then  $\calF_\la$ for $\la\geq 1$ is defined at
least over $(-\delta/\sqrt{\la},\delta/\sqrt{\la})$.
\end{conjecture}
}

\begin{proof}   Combining \eqref{free cumulants}
with $R_{aX+b}(z)=b+a R_X(az)$, we see that the free cumulants of
$\nu_{\la}$ are $c_1(\nu_{\la})=c_1(\nu)=m_0$ and for $n\geq 1$
\begin{equation*}
c_{n+1}(\nu_{\la})=\frac{1}{\la^{n}}c_{n+1}(\nu_{}) =\left.\frac{1}{
n!}\frac{d^{n-1}}{d
x^{n-1}}\left(\frac{V(x)}{\la}\right)^n\right|_{x=m_0}.
\end{equation*}
Theorem \ref{V2c} implies that $V/\la$ is the variance function of
the free exponential  family generated by $\nu_{\la}$.

If $\nu_{1/n}$ exists for all $n\in\NN$, then the first part of the
proposition together with uniqueness theorem (Theorem \ref{T1})
implies that
 $\nu_{}=(D_{n}(\nu_{1/n}))^{\boxplus n}$, proving $\boxplus$-infinite divisibility.
\end{proof}

\section{Marchenko-Pastur Approximation}\label{Sect FPA}
Let $$\omega_{a,\sigma}(dx)=\frac{\sqrt{4\sigma^2-(x-a)^2}}{2\pi
\sigma^2}1_{|x-a|<2\sigma}dx$$ denote the semicircle law of mean $a$
and variance $\sigma^2$. Up to affine transformations, this is the
free Meixner law which appears in Theorem \ref{T free quadr} as the
law which  generates the free exponential  family  $\calF_a(V)$ with
the variance function $V \equiv \sigma^2$.

{\tolerance=500
Following the analogy with natural exponential families, family
$\calF_0(\sigma^2)$
 can be thought as a free exponential analog of the normal family.
Somewhat surprisingly, this   family does not contain all semicircle laws,
 but instead it
contains affine transformations of the (absolutely continuous) Marchenko-Pastur laws.}
\begin{example}[Semi-circle free exponential family]\label{Example: semicircle}
For $\la>0$, let
$$\pi_{m,\la}(dx)=\frac{\sqrt{\la} \sqrt{4-\la
x^2}}{2\pi (1+\la m (m-x))}1_{x^2\leq 4/\la} dx.$$
Function $V(m)\equiv1/\la$ is the  variance function of
the free exponential family \begin{equation}\label{pi}
\calF_0({1/\la})=\left\{\pi_{m,\la}(dx): |m|< 1/\sqrt{\la}\right\}
\end{equation}
with the generating measure $\nu(dx)= \omega_{0,1/\sqrt{\la}}(dx)$.

To verify that the expression integrates to $1$ for  $m\ne 0$,   we use the explicit form of the  density  \cite[(3.3.2)]{Hiai-Petz00}  to note that
$\pi_{m,\la}=\calL(m+1/(\la m)- mX)$ is the law of the affine transformation of a free Poisson (Marchenko-Pastur)
 random variable $X$ with parameter $1/(\la m^2)$.
From the properties of Marchenko-Pastur law we see that $\int
\pi_{m,\la}(dx)=1$ iff $m^2\leq 1/\la$, so for large $\la$ the
interval $(A,B)\subset(-1/\sqrt{\la},1/\sqrt{\la})$ in \eqref{F(V)}
cannot be chosen  independently of $\la$.

We remark that  Biane \cite{Biane-97}  analyzes $f\mapsto
g(m):=\int f(x)\pi_{m,\la}(dx)$  as a mapping of the
appropriate Hilbert spaces for complex  $m$.
\end{example}

We have the following analogue  of \cite[Theorem 3.4]{Jorg}.
\begin{theorem}[Marchenko-Pastur approximation]\label{T-CLT}
Suppose   the  variance function $V$ of a free exponential  family
$\mathcal{F}_{m_0}(V)$ is analytic and strictly positive  in a neighborhood of $m_0$. Then
there is $\delta>0$ such that if $\mathcal{L}(Y_\la)\in
\mathcal{F}_{m_0}(V/\la)$ has mean
 $E(Y_\la)=m_0+m/\sqrt{\la}$ with $|m|<\delta$, then
$\sqrt{\la}(Y_\la-m_0) \xrightarrow{\calD} \pi_{m,1/V(m_0)}$ as
$\la\to\infty$.
\end{theorem}

To prove Theorem \ref{T-CLT} we will use  the following analogue of Mora's
Theorem, see \cite[Theorem 2.12]{Jorg}, or \cite[Theorem
2.6]{Let:Mor}.
\begin{proposition}\label{Thm free Mora} \tolerance=1000
Suppose $V_n$ is a family of analytic functions which are variance functions of  a sequence of
free exponential  families  $\left\{\mathcal{F}_{m_0}(V_n):n\geq 1\right\}$.
If $V_n\to V$ uniformly in a (complex) neighborhood of  $m_0\in\RR$, and $V(m_0)>0$, then there is $\delta>0$ such that $V$
is a  variance function of a free exponential  family
$\mathcal{F}_{m_0}(V)$ parameterized by the mean
$m\in(m_0-\delta,m_0+\delta)$. Moreover, if a sequence of measures
$\mu_n\in\mathcal{F}_{m_0}(V_n)$ is such that $m_1=\int
x\mu_n(dx)\in(m_0-\delta,m_0+\delta)$ does not depend on $n$, then
$\mu_n\xrightarrow{\calD} \mu$ where $\mu\in \mathcal{F}_{m_0}(V)$
has the same mean $\int x\mu(dx)=m_1$.
\end{proposition}
\begin{proof} Let $\nu_n$ be the generating measure for $\mathcal{F}_{m_0}(V_n)$.
Since $V_n(z)\to V(z)$  uniformly in a neighborhood of $m_0$,
from  \eqref{Cauchy}
we see that the cumulants $c_{k+1}(\nu_n)$ converge as $n\to\infty$
and $\sup_n |c_{k+1}(\nu_n)|\leq M^k$ for some $M<\infty$. Therefore
 the $R$-transforms of $\nu_n$  converge to the $R$-transform of a
compactly supported measure $\nu$. Thus
$\nu_n\xrightarrow{\calD}\nu$, and the supports of $\nu_n$ are
uniformly bounded in $n$, i.e.,
 $\mbox{supp}(\nu_n)\subset[-A,A]$ for some $0<A<\infty$.
By decreasing the value of  $\delta$ we can also ensure that the
densities in \eqref{Q_m} are  bounded as functions of $x\in[-A,A]$
uniformly in $n$. So the integrals converge, and $\nu$ indeed generates a free exponential family with variance $V$ in a neighborhood of $m_0$.

Suppose now $\mu_n\in\mathcal{F}_{m_0}(V_n)$  and
$\mu\in \mathcal{F}_{m_0}(V)$ have the same
mean $m$ for some $|m-m_0|<\delta$ small enough.  Since the
densities  in \eqref{Q_m}  are  bounded  by some constant $C$ for
all $x\in[-A,A]$, $n\in\NN$, for a bounded continuous function $g$
we have
\begin{multline*}
\int g(x)\mu_n(dx)=\int g(x)\frac{V_n(m)}{V_n(m)+ (m-m_0)(m-x)}\nu_n(dx)\\
=
\int g(x)\frac{V(m)}{V(m)+ (m-m_0)(m-x)}\nu_n(dx)+\eps_n\\
\to\int g(x)\frac{V(m)}{V(m)+ (m-m_0)(m-x)}\nu(dx)= \int
g(x)\mu(dx)
\end{multline*}
as $n\to \infty$, where
\begin{multline*}
|\eps_n|=\Big|\int g(x)\big(\frac{V_n(m)}{V_n(m)+ (m-m_0)(m-x)}\\
-\frac{V(m)}{V(m)+ (m-m_0)(m-x)}\big)\nu_n(dx)\Big|\\
\leq C^2 | m-m_0|\sup_{x\in [-A,A]} |g(x)(m-x)|\frac{|V_n(m)-V(m)|}{V_n(m)V(m)}\to 0.
\end{multline*}

\end{proof}
\begin{proof}[Proof of Theorem \ref{T-CLT}]
Without loss of generality, we assume $m_0=0$.  Suppose
$\mu\in\calF(V)$. A change of variable shows that
$D_a(\mu)(dx)\in\calF(V(a m)/a^2)$,
 with generating measure $\nu_a=D_a(\nu)$.
Since
$$\calL(Y_\la)\in \calF\left(\frac{V(m)}{\la}\right),$$
 this shows that
  $$\calL(\sqrt{\la}Y_\la)\in\calF\left(V\left(\frac{m}{\sqrt{\la}}\right)\right).$$
\tolerance=700 We now use Proposition \ref{Thm free Mora} to the
sequence of variance functions
$V_\la(m)=V\left(\frac{m}{\sqrt{\la}}\right)\to V(0)=\sigma^2$ as
$\la\to\infty$. From Proposition \ref{Thm free Mora} we deduce that
there is $0<\delta<\sigma^2$ such that if $|m|<\delta$ and
 $E(\sqrt{\la}Y_\la)=m$, then   $\calL(\sqrt{\la}Y_\la)\xrightarrow{\calD} \pi_{m,1/\sigma^2}\in\calF_0(\sigma^2)$, see \eqref{pi}.
\end{proof}
By Example \ref{Example: semicircle}, if $0<|m|\leq \sigma$, then up
to affine transformation $\pi_{m,1/\sigma^2}$ is a Marchenko-Pastur
law. Thus in this case Theorem \ref{T-CLT} gives a Marchenko-Pastur
approximation to $\calL(Y_\la)$.\tolerance=500

Of course, every compactly supported mean-zero measure $\nu$ is an element of the
Cauchy-Stieltjes family that it generates. Since
$\pi_{0,1/\sigma^2}=\omega_{0,\sigma}$ is the semicircle law,
combining Proposition \ref{C.4.3} with Theorem \ref{T-CLT} we get
the following Free Central Limit Theorem; see \cite{Bozejko-75,Voiculescu86}.
\begin{corollary}
If a probability measure $\nu$ is compactly supported and centered, then
with $\sigma^2=\int x^2\nu(dx)$ we have
$$
D_{\sqrt{n}}(\nu^{\boxplus n})\xrightarrow{\calD}\omega_{0,\sigma}\;.
$$
\end{corollary}

\subsection*{\bf Acknowledgements} The   author  thanks J. Weso\l owski
for the copy of \cite{Wesolowski90} and for several helpful discussions.
He also thanks the referee for an insightful report which helped  to clarify the statement of Theorem \ref{T1} and led  to significant improvement of the paper.


\begin{thebibliography}{10}

\bibitem{Anshelevich01}
M.~Anshelevich.
\newblock Free martingale polynomials.
\newblock {\em J. Funct. Anal.}, 201:228--261, 2003.
\newblock arXiv:math.CO/0112194.

\bibitem{Anshelevich04}
M.~Anshelevich.
\newblock Orthogonal polynomials with a resolvent-type generating function.
\newblock {\em Trans. Amer. Math. Soc.}, 360(8):4125--4143, 2008.
\newblock arXiv:math.CO/0410482.

\bibitem{Barndorff-Nielsen-Thorbj02a}
O.~E. Barndorff-Nielsen and S.~Thorbj{\o}rnsen.
\newblock Self-decomposability and {L}\'evy processes in free probability.
\newblock {\em Bernoulli}, 8(3):323--366, 2002.

\bibitem{Ben_Arous:2009}
G.~Ben~Arous and V.~Kargin.
\newblock Free point processes and free extreme values.
\newblock {\em Probability Theory and Related Fields}, (to appear),
  arXiv:0903.2672, 2009.

\bibitem{Ben-Arous:2006}
G.~Ben~Arous and D.~V. Voiculescu.
\newblock Free extreme values.
\newblock {\em Ann. Probab.}, 34(5):2037--2059, 2006.

\bibitem{Benaych-Georges04}
F.~Benaych-Georges.
\newblock Taylor expansions of ${R}$-transforms, application to supports and
  moments.
\newblock {\em Indiana Univ. Math. J.}, 55:465--482, 2006.

\bibitem{Biane-97}
P.~Biane.
\newblock Segal-{B}argmann transform, functional calculus on matrix spaces and
  the theory of semi-circular and circular systems.
\newblock {\em J. Funct. Anal.}, 144(1):232--286, 1997.

\bibitem{Bozejko-75}
M.~Bo{\.z}ejko.
\newblock On {$\Lambda (p)$} sets with minimal constant in discrete
  noncommutative groups.
\newblock {\em Proc. Amer. Math. Soc.}, 51:407--412, 1975.

\bibitem{Bozejko-Bryc-04}
M.~Bo\.zejko and W.~Bryc.
\newblock On a class of free {L}\'evy laws related to a regression problem.
\newblock {\em J. Funct. Anal.}, 236:59--77, 2006.
\newblock arxiv.org/abs/math.OA/0410601.

\bibitem{Bryc-Ismail-05}
W.~Bryc and M.~Ismail.
\newblock Approximation operators, exponential, and $q$-exponential families.
\newblock Preprint. arxiv.org/abs/math.ST/0512224, 2005.

\bibitem{Bryc-Wesolowski-03}
W.~Bryc and J.~Weso{\l}owski.
\newblock Conditional moments of $q$-{M}eixner processes.
\newblock {\em Probab. Theory Related Fields}, 131:415--441, 2005.
\newblock arxiv.org/abs/math.PR/0403016.

\bibitem{Diaconis-08}
P.~Diaconis, K.~Khare, and L.~Saloff-Coste.
\newblock Gibbs sampling, exponential families and orthogonal polynomials.
\newblock {\em Statistical Science}, 23:151--178, 2008.

\bibitem{Hiai-Petz00}
F.~Hiai and D.~Petz.
\newblock {\em The semicircle law, free random variables and entropy},
  volume~77 of {\em Mathematical Surveys and Monographs}.
\newblock American Mathematical Society, Providence, RI, 2000.

\bibitem{Ism:May}
M.~E.~H. Ismail and C.~P. May.
\newblock On a family of approximation operators.
\newblock {\em J. Math. Anal. Appl.}, 63(2):446--462, 1978.

\bibitem{Jorg}
B.~J{\o}rgensen.
\newblock {\em The theory of dispersion models}, volume~76 of {\em Monographs
  on Statistics and Applied Probability}.
\newblock Chapman \& Hall, London, 1997.

\bibitem{Let}
G.~Letac.
\newblock {\em Lectures on natural exponential families and their variance
  functions}, volume~50 of {\em Monograf\'\i as de Matem\'atica [Mathematical
  Monographs]}.
\newblock Instituto de Matem\'atica Pura e Aplicada (IMPA), Rio de Janeiro,
  1992.

\bibitem{Let:Mor}
G.~Letac and M.~Mora.
\newblock Natural real exponential families with cubic variance functions.
\newblock {\em Ann. Statist.}, 18(1):1--37, 1990.

\bibitem{Mor}
C.~N. Morris.
\newblock Natural exponential families with quadratic variance functions.
\newblock {\em Ann. Statist.}, 10(1):65--80, 1982.

\bibitem{Nica-Speicher}
A.~Nica and R.~Speicher.
\newblock On the multiplication of free {$N$}-tuples of noncommutative random
  variables.
\newblock {\em Amer. J. Math.}, 118(4):799--837, 1996.

\bibitem{Pol:Sze}
G.~P{\'o}lya and G.~Szeg{\H{o}}.
\newblock {\em Problems and theorems in analysis. {I}}, volume 193 of {\em
  Grundlehren der Mathematischen Wissenschaften [Fundamental Principles of
  Mathematical Sciences]}.
\newblock Springer-Verlag, Berlin, 1978.

\bibitem{Saitoh-Yoshida01}
N.~Saitoh and H.~Yoshida.
\newblock The infinite divisibility and orthogonal polynomials with a constant
  recursion formula in free probability theory.
\newblock {\em Probab. Math. Statist.}, 21(1):159--170, 2001.

\bibitem{Speicher-97}
R.~Speicher.
\newblock Free probability theory and non-crossing partitions.
\newblock {\em S\'em. Lothar. Combin.}, 39:Art.\ B39c, 38 pp.\ (electronic),
  1997.

\bibitem{Voiculescu86}
D.~Voiculescu.
\newblock Addition of certain noncommuting random variables.
\newblock {\em J. Funct. Anal.}, 66(3):323--346, 1986.

\bibitem{Voiculescu00}
D.~Voiculescu.
\newblock Lectures on free probability theory.
\newblock In {\em Lectures on probability theory and statistics (Saint-Flour,
  1998)}, volume 1738 of {\em Lecture Notes in Math.}, pages 279--349.
  Springer, Berlin, 2000.

\bibitem{Wesolowski90}
J.~Weso{\l}owski.
\newblock Kernel families.
\newblock Unpublished manuscript, 1999.

\end{thebibliography}
\end{document}